\input amssym.def
\input amssym.tex
\input pictex
\magnification=\magstep1
\hfuzz=5pt
\def\cT{{\cal T}}
\def\V{{\bf V}} 
\font\titl=cmbx10 scaled\magstep2      
\font\auth=cmsl10 scaled\magstep1      
\font\bfootfont=cmr8                   
\font\fit=cmti8                        
%
%
\font\tenmsa=msam10
\font\sevenmsa=msam7
\font\fivemsa=msam5
\newfam\msafam
\textfont\msafam=\tenmsa
\scriptfont\msafam=\sevenmsa
\scriptscriptfont\msafam=\fivemsa

\font\tenmsb=msbm10
\font\sevenmsb=msbm7
\font\fivemsb=msbm5
\newfam\msbfam
\textfont\msbfam=\tenmsb
\scriptfont\msbfam=\sevenmsb
\scriptscriptfont\msbfam=\fivemsb

\mathchardef\subsetneqq="2924 
\mathchardef\supsetneqq="2925
%
%
\def\ABC{{\bf 1}} 
\def\BGSS{{\bf 2}} 
\def\B{{\bf 3}}
\def\Can{{\bf 4}}
\def\C{{\bf 5}} 
\def\CMSZa{{\bf 6}} 
\def\Ch{{\bf 7}}
\def\CD{{\bf 8}}
\def\Co{{\bf 9}}
\def\D{{\bf 10}}
\def\Ep{{\bf 11}}
\def\GS{{\bf 12}}
\def\GSrat{{\bf 13}}
\def\Gr{{\bf 14}}
\def\Mou{{\bf 15}}
\def\Pa{{\bf 16}}
\def\Pap{{\bf 17}}
\def\Ro{{\bf 18}}
\def\T{{\bf 19}} 

{\titl \multiply\baselineskip by 130\divide\baselineskip by 100
\leftline{Hyperbolic buildings, affine buildings}
\medskip
\leftline{and automatic groups}}
\vskip 1em

\leftline{\auth Donald I. Cartwright\footnote{$^1$}
  {\bfootfont School of Mathematics and Statistics, University of Sydney, 
   N.S.W. 2006, Australia.}
   and Michael Shapiro\footnote{$^2$}{\bfootfont Department of
   Mathematics, City College of New York, New York, NY 10031, USA \hfil\break
   {\fit Key words and phrases.} Automatic groups, affine buildings,
   hyperbolic groups.}} \vskip2pt \vskip 1.5em

{\narrower\noindent {\bf Abstract.} We see that a building whose
Coxeter group is hyperbolic is itself hyperbolic.  Thus any finitely
generated group acting co-compactly on such a building is hyperbolic,
hence automatic.  We turn our attention to affine buildings and
consider a group $\Gamma$ which acts simply transitively and in a
``type-rotating" way on the vertices of a locally finite thick
building of type~$\widetilde{A_n}$. We show that $\Gamma$ is
biautomatic, using a presentation of~$\Gamma$ and unique normal form
for each element of~$\Gamma$, as described in~[\C].  \vskip 1.5em}

\noindent{\bf\S1. Introduction.}  

\medskip\noindent Two standard references for the theory of automatic
groups are the book~[\Ep] and the paper~[\BGSS], which both contain
numerous examples of automatic groups.  Perhaps the ``canonical''
class of automatic groups is the class of word hyperbolic groups of
Gromov {\it et al.}~[\Gr]. As we shall see in Section 3, any building
whose underlying Coxeter group is word hyperbolic is itself hyperbolic
in the word metric.  In particular any finitely generated group acting
co-compactly with finite stabilizers on such a building is word
hyperbolic and thus automatic.  Thus, from the viewpoint of automatic
groups, it is natural to look next at actions on affine buildings.  

The first result in this direction is provided by Gersten and
Short~[\GS] who prove that a finitely generated torsion free group
that acts co-compactly discretely by isometries on a Euclidean
building of dimension 2 is automatic.  It has often seemed likely that
the restriction on dimension could be lifted. In this
paper, we prove this for finitely generated groups which act simply
transitively (and in a type-rotating way) on the vertices of a thick
building of type~$\widetilde{A_n}$. We call such groups {\it
$\widetilde{A_n}$-groups\/}.  In fact, we will show more.  We shall
see that these groups are biautomatic.  The structure in question is a
symmetric automatic structure~[\Ep].  (The term ``fully automatic'' is
used in [\Ch] and early versions of [\Ep].)  In particular this
implies that  $\widetilde{A_n}$-groups have solvable conjugacy
problem~[\GSrat].  Further, the structure consists of geodesics, and
thus these groups have rational growth functions.  Examples of
finitely generated $\widetilde{A_n}$-groups are known in dimensions
$n=2~[\CMSZa]$, $n=3,4~[\C]$ and $n=5$~[Cartwright, unpublished], and
for $n=2,3,4$ and  for any prime power $q$ there are examples of
$\widetilde{A_n}$-groups which are arithmetic lattices in
$PGL(n+1,{\bf F}_q((X)))$.   It remains an open problem whether
$\widetilde{A_n}$-groups exist for every~$n$.

The paper is organized as follows.  Section~2 defines hyperbolic groups
and reviews some basic background information, most of which can be found
in [\ABC].  In Section 3 we show that a building is hyperbolic (in an
appropriate sense) if and only if its underlying Coxeter group is
hyperbolic.   Section 4 defines $\widetilde{A_n}$-groups and reviews
the necessary background from [\C].  Section 5 gives the proof that
finitely generated $\widetilde{A_n}$-groups are bi-automatic.
Sections 4 and 5 require no knowledge of buildings, Section 3 uses
elementary results which can be found in [\B, IV.3].

\medskip\noindent{\bf\S2.  Hyperbolic groups}\medskip

\noindent 
We say that a metric space $(X,d)$ is a {\it geodesic metric space} if
for every $x,y\in X$, there is a path from $x$ to $y$ which realizes
their distance.  Such a path is called a {\it geodesic}.  Following
[\Gr] we say that a geodesic metric space is {\it
$\delta$-hyperbolic} if whenever $P$ is a point on side $\alpha$ of a
geodesic triangle with sides $\alpha$, $\beta$, and $\gamma$, there
is a point $Q$ on $\beta \cup \gamma$ so that $d(P,Q) \le \delta$.
We say that $(X,d)$ is {\it hyperbolic} if it is $\delta$-hyperbolic
for some $\delta$.

Now given any connected graph $\Gamma$ there is a natural metric on
$\Gamma$.  Take each edge of $\Gamma$ to be isometric to the
unit interval and take the path metric which this induces on
$\Gamma$.  Further, given any finitely generated group $G$, the
choice of a finite generating set ${\cal G}=\{a_1,\ldots,a_k\}$ turns
$G$ into a  directed labelled connected graph $\Gamma=\Gamma_{\cal G}$.  The
vertices of $\Gamma$ are the elements of $G$ and the edges of
$\Gamma$ are $\{(g,ga) \mid g \in G, \ a \in \cal G\}$.  We direct
the edge $(g,ga)$ from $g$ to $ga$ and label it with $a$.  We assume
that $\cal G$ is closed under inverses and identify $(g,ga)$ with the
inverse of $(ga,g)$.  $\Gamma$ is called the {\it Cayley graph} of
$G$ with respect to $\cal G$.

We say that $G$ is {\it hyperbolic} if $\Gamma=\Gamma_{\cal G}$ is
hyperbolic.  While this appears to depend on $\cal G$, in fact only
the particular value of $\delta$ depends on $\cal G$.  

We will want a standard fact about hyperbolic metric spaces.  Given a
path $\sigma$ in $X$ and $0 < \lambda \le 1$ and $0 \le \epsilon$, we
will say that $\sigma$ is a {\it $(\lambda,\epsilon)$-quasigeodesic}
if for every decomposition $\sigma=\alpha\beta\gamma$, the endpoints
of $\beta$ are separated by at least $\lambda\ell(\beta)-\epsilon$.
(Here $\ell(\beta)$ denotes the length of $\beta$.)  If $X$ is a
$\delta$-hyperbolic metric space then there is
$N=N(\delta,\lambda,\epsilon)$ so that if $\sigma$ is a
$(\lambda,\epsilon)$-quasigeodesic and $\tau$ is a geodesic with the
same endpoints, then $\sigma$ and $\tau$ each lie in the
$N$-neighborhood of each other [\ABC, 3.3].  Geodesics are simply
$(1,0)$-quasigeodesics, and thus there is $N=N(\delta,1,0)$ so that
all geodesics joining common endpoints live in a $N$-neighborhood of
each other.  We call a pair of geodesics with common endpoints a
{\it bigon}. 

 From this one can construct a proof that any  geodesic metric space
quasi-isometric to a hyperbolic space is itself hyperbolic, and that
in particular, hyperbolicity of a group is independent of generating
set.  Now it is a standard result that if a finitely generated group
$G$ acts co-compactly by isometries and with finite stabilizers on a
geodesic metric space $(X,g)$, then every Cayley graph of $G$ is
quasi-isometric to $X$. (See, for example, [\Can].)  In particular,
when $X$ is hyperbolic, so is $G$.

We have seen that the definition of a hyperbolic metric space requires
checking that geodesic triangles are ``thin'', that is, that no side
of a triangle is ever far from the union of the other two sides.
Papasoglu [\Pa], [\Pap], has shown that in graphs it is only necessary
to check bigons.   Note that the endpoints of geodesics (and hence the
endpoints of bigons) need not be vertices of $\Gamma$.

\medskip\noindent{\bf Theorem.} (Papasoglu)
{\sl Suppose that $\Gamma$ is a graph and there is a constant $K$ so
that if $\sigma$, $\sigma'$ is a bigon, then $\sigma$ and $\sigma'$
each lie in a $K$-neighborhood of each other.  Then $\Gamma$ is
hyperbolic. }\medskip

Notice that the hypothesis is equivalent to the {\it a priori}
stronger hypothesis that there exists $K'$ so that if $\sigma$ and
$\sigma'$ form a bigon then for all $t$, $d(\sigma(t), \sigma'(t)) \le
K'$.  For suppose that $\sigma(t)$ is within $K$ of $\sigma'(t')$.
Then the fact that these are geodesics emanating from a common point
allows us to use the triangle inequality to see that $|t-t'|\le K$.
Hence, taking $K'=2K$, we have $d(\sigma(t), \sigma'(t)) \le
K'$.

\medskip\noindent{\bf\S3.  Hyperbolic buildings}\medskip

\noindent In his doctoral thesis [\Mou], Moussong constructs actions
of Coxeter groups on  non-positively curved geodesic metric spaces.
An account of this metric can be found in [\Co].  The metric spaces in
question are locally Euclidean or locally hyperbolic complexes, and
can be made negatively curved if and only if the Coxeter group in
question is word hyperbolic.  The actions are co-compact, by
isometries and with finite stabilizers.  As a scholium of his
construction, one knows exactly which Coxeter groups are word
hyperbolic.

\medskip\noindent{\bf Theorem} (Moussong) {\sl Let $(W,S)$ be a
Coxeter system. Then the following are equivalent:

\item{\it 1.} $W$ is word hyperbolic.

\item{\it 2.} $W$ has no ${\bf Z}\times {\bf Z}$ subgroup.

\item{\it 3.} $(W,S)$ does not contain a affine sub-Coxeter system of
rank $\ge 3$, and does not contain a pair of disjoint commuting
sub-Coxeter systems whose groups are both infinite.} \smallskip

Charney and Davis [\CD] have pointed out that using Moussong's metric of
non-positive curvature, one can give a building a metric of
non-positive curvature, and that the metric on the building is
negatively curved if and only if the Coxeter group is word hyperbolic.
(Construction of the metric on the building can be done along the
lines of [\B, VI.3].)

We give a similar characterization in terms of graphs.  Given a
building $\Delta$, there is a metric on the set of chambers of
$\Delta$, and we will want a path metric space which reflects this
metric.  To do this we let $\Delta'$ be the graph dual to $\Delta$.
That is to say, the vertices of $\Delta'$ are the barycenters of the
chambers of $\Delta$.  Two such vertices are connected by an edge when
they lie in chambers with a common face.  As usual, $\Delta'$ is
metrized considering each edge as isometric to the unit interval.
Non-stuttering galleries of $\Delta$ correspond to edgepaths in
$\Delta'$. The decomposition of $\Delta$ into apartments induces a
decomposition of $\Delta'$ into apartments which are isometric as
labelled graphs to the Cayley graph of $(W,S)$, the Coxeter system of
$\Delta$.

\medskip\noindent{\bf Theorem 1.} {\sl Suppose $\Delta$ is a building whose
apartments are the Coxeter complex of a word hyperbolic Coxeter group.
Then $\Delta'$ is hyperbolic.}

\medskip\noindent{\bf Remark.}  The converse is also true.  That is, if
$\Delta'$ is hyperbolic, the associated Coxeter group is word
hyperbolic.  This follows immediately from the fact that the embedding
of the Cayley graph into $\Delta'$ is an isometry.

\medskip\noindent{\bf Corollary.} {\sl
Suppose $\Delta$ is a building whose apartments are the Coxeter complex
of a word hyperbolic Coxeter group and that $G$ is a finitely
generated group which acts simplicially, co-compactly with finite
stabilizers on $\Delta$.  Then $G$ is word hyperbolic.}

\medskip\noindent{\bf Proof.}
$G$ acts on $\Delta$, and the natural embedding of $\Delta'$ into
$\Delta$ is equivariant with respect to this action.  This induces an
action of $G$ on $\Delta'$.  Since there are finitely many $G$-orbits
of chambers of $\Delta$ and finitely many $G$-orbits of codimension 1
faces, the induced action on $\Delta'$ is co-compact.  Now $G$ carries
edges of $\Delta'$ to edges of $\Delta'$ and thus acts by isometries
of the graph metric.  Since the action of $G$ has finite stabilizers,
the (setwise) stabilizer of each chamber and co-dimension 1 face of
$\Delta$ is finite.  These are the stabilizers of the vertices and
edges of $\Delta'$.  Hence, as we have outlined in Section 2, $G$ is
quasi-isometric to $\Delta'$ and thus word hyperbolic.

\medskip\noindent{\bf Proof of Theorem 1.}
To prove the Theorem, we consider a bigon $\sigma$, $\sigma'$.  It
suffices to show that there is $K$ so that if $C$ is any point of
$\sigma$, then $C$ is within $K$ of $\sigma'$.  We distinguish three
cases depending on whether both, one, or neither  of the endpoints of
this bigon are vertices.

Case 1:  Both ends of the bigon are vertices.  Now  given any two
chambers of $\Delta$, there is an apartment $\Sigma$ containing them
both. Any geodesic gallery in $\Delta$ connecting these two chambers
lies in $\Sigma$.  (See [\B, p.~88].)  It now follows that $\sigma$
and $\sigma'$ lie in a common apartment of $\Delta'$, and by the
hyperbolicity of the underlying Coxeter group, we are done.

Case 2: The geodesics $\sigma$ and $\sigma'$ begin at a vertex, but do
not end at a vertex.  We let $x$ and $x'$ be the beginning and
endpoints of our bigon.  We let $y$ and $z$ be the last vertices of
$\sigma$ and $\sigma'$ respectively, and let $\tau$ and $\tau'$ be the
initial segments of $\sigma$ and $\sigma'$ ending at $y$ and $z$
respectively.   It now follows that $\ell(\tau)=\ell(\tau')$ and that
$x'$ is the midpoint of an edge $e$.  If $y=z$, $\tau$ and $\tau'$ lie
in a common apartment, and we are done by Case~1. Thus we may assume
$y\ne z$.  Now $\tau$ and $\tau'$ cannot lie in a common apartment.
For if this were so, $\tau e \tau'^{-1}$ would label a relator of odd
length in the underlying Coxeter group, and this is impossible. 
$$
\beginpicture
\setcoordinatesystem units <1mm,1mm>
\setplotarea x from -20 to 50, y from -5 to 20
\setlinear
\plot -8.66 0 -8.66 10 0 5 -8.66 0 / 
\plot -2.89 10 -5.77 5 -2.89 0 /
\plot 40 5 45 13.66 50 5 45 -3.66 40 5 50 5 /
\plot 40 -0.88 45 2.22 45 7.78 40 10.88 /
\put {$\scriptstyle\bullet$} at 45 7.78
\put {$\scriptstyle\bullet$} at 45 2.22
\put {$\scriptstyle\bullet$} at -5.77 5
\put {$\scriptstyle\bullet$} at  20 15
\put {$\scriptstyle\bullet$} at  40 -0.88
\put {$\scriptstyle\bullet$} at  40 10.88
\put {$\scriptstyle\bullet$} at  -2.89 0
\put {$\scriptstyle\bullet$} at  -2.89 10
\put {$x$} at -7 6
\put {$y$} at 46.5 7.78
\put {$z$} at 46.5 2.22  
\put {$C$} at 21 17.5  
\put {$\sigma$} at 11 15.5  
\put {$\sigma'$} at 11 -2.5 
\put {$x'$} at 52 8
\put {Case 2} at 20 -10 
\arrow <6pt> [.15,.6] from 50 7 to 45 5 
\setdashes
\setquadratic 
\plot -2.89 10 10 14 20 15 30 14 40 10.88 /
\plot -2.89  0 10 -4 20 -5 30 -4 40 -0.88 /
\endpicture 
$$
Let $\Sigma$ be an apartment containing $\sigma$. We consider
$C=\sigma(t)$ and may suppose this  is a vertex lying  on $\tau$.  We
let $\rho=\rho_{\Sigma,C}$ be the retraction onto $\Sigma$ centered at
$C$.  (See, for example  [\B, IV.3].)  Now $d(y,\rho(z))\le 1$ and
since $\rho$ does not increase distance, we have $$n-1 \le d(x,
\rho(z)) \le n,$$ where $n=\ell(\tau)$.  Now consider the path
$\rho(\tau')$.  This may not be an edge path, as some edge of $\tau'$
may be folded by $\rho$. However its length as a path is still $n$.
We thus have a path of  length $n$ whose endpoints lie at distance at
least $n-1$.  It follows that $\rho(\tau')$ is a
$(1,1)$-quasigeodesic.  Since $\rho(\tau')$ lies in $\Sigma$,  it now
follows from the hyperbolicity of $\Sigma$ that $\rho(\tau')$ lies
close to $\tau$.  We can thus find $C'=\tau'(t')$ so that $\rho(C')$
lies close to $C$.  Since $\rho=\rho_{\Sigma,C}$ preserves distance
from $C$, $C'$ lies close to $C$ and we are done.

Case 3:  Neither end of the bigon is a vertex.  In this case we let
$a$ and $b$, (respectively $a'$ and $b'$) be the first and last
vertices of $\sigma$ (respectively $\sigma'$), and let $\tau$
(respectively $\tau'$) be the segment of $\sigma$ (respectively
$\sigma'$) connecting these first and last vertices.  If $a=a'$ or
$b=b'$ or both we are reduced to previous cases, so we can assume $a
\ne a'$ and $b \ne b'$.  We let $e$ be the edge from $a$ to $a'$ and
$e'$ be the edge from $b$ to $b'$.  It is easy to check that
$\ell(\tau)=\ell(\tau')$.  We take $n=\ell(\tau)$.  Thus $d(a,b')$ is
either $n-1$, $n$ or $n+1$.  We let $\mu$ be a geodesic from $a$ to
$b'$. 
$$
\beginpicture
\setcoordinatesystem units <1mm,1mm>
\setplotarea x from -20 to 50, y from -5 to 20
\setlinear
\plot 40 5 45 13.66 50 5 45 -3.66 40 5 50 5 /
\plot 40 -0.88 45 2.22 45 7.78 40 10.88 / 
\plot -10 5 -5 13.66 0 5 -5 -3.66 -10 5 0 5 / 
\plot 0 -0.88 -5 2.22 -5 7.78 0 10.88 / 
\put {$\scriptstyle\bullet$} at 45 7.78
\put {$\scriptstyle\bullet$} at 45 2.22
\put {$\scriptstyle\bullet$} at 40 10.88
\put {$\scriptstyle\bullet$} at 40 -0.88
\put {$\scriptstyle\bullet$} at  20 15
\put {$\scriptstyle\bullet$} at  0 -0.88
\put {$\scriptstyle\bullet$} at  0 10.88
\put {$\scriptstyle\bullet$} at  -5 7.78
\put {$\scriptstyle\bullet$} at  -5 2.22
\put {$b$}  at 46.5 8.25
\put {$b'$} at 47 2.22  
\put {$a$}  at -6.5 7.78
\put {$a'$} at -6.7 2.7  
\put {$e'$} at 51.5 8
\put {$e$} at -11 8
\put {$C$} at 21 17.5  
\put {$\sigma$} at 11 15.5  
\put {$\sigma'$} at 11 -2.5   
\put {$\mu$} at 21 6.5  
\put {Case 3} at 20 -10 
\arrow <6pt> [.15,.6] from 50 7 to 45 5.5 
\arrow <6pt> [.15,.6] from -10 7 to -5 5.5 
\setdashes
\setquadratic 
\plot 0 10.88 10 14 20 15 30 14 40 10.88 /
\plot 0 -0.88 10 -4 20 -5 30 -4 40 -0.88 /
\plot -5 7.78 10 6.5 20 5 30 4 45 2.22 /
\endpicture 
$$
If $d(a,b')=n-1$ then $\mu e'^{-1}$ and $\tau$ form a bigon
whose ends are the vertices $a$ and $b$.  Likewise $e^{-1}\mu$ and
$\tau'$ form a bigon whose ends are the vertices $a'$ and $b'$.
Applying Case 1 twice takes care of this situation.

If $d(a,b')= n$, we let $y$ and $y'$ be the midpoints of $e$ and
$e'$.  We then have a bigon whose ends are $a$ and $y'$ whose sides
consist of $\tau$ and $\mu$, each with a half of $e'$ appended.
Similarly, we have a bigon with ends $b'$ and $y$ and sides
$\mu^{-1}$ and $\tau'^{-1}$ with the halves of $e$ appended.  Now we
can apply Case 2 twice.

Finally, if $d(a,b')=n+1$ then $\mu$ and $\tau e'$ form a bigon as do
$\mu$ and $e\tau'$, and once again we can apply Case 1 twice.

\medskip\noindent{\bf\S4. Review of $\widetilde A_n$-groups.} 
\medskip 

\noindent $\widetilde{A_n}$-groups were introduced for general $n \ge
2$ in~[\C], after earlier work on the case $n=2$ in~[\CMSZa].  (One
dimensional buildings are trees, and thus hyperbolic.) Recall that a
building is a {\it labellable} complex, that is, that each vertex~$v$
of a building~$\Delta$ of type~$\widetilde{A_n}$ has a {\it type\/}
$\tau(v)\in\{0,1,\ldots,n\}$, with each chamber having one vertex of
each type.  If $g$ is an automorphism of~$\Delta$,
and if there is an integer~$c$ such that
$\tau(gv)=\tau(v)+c$~(mod~$n+1$) for each vertex $v$, then $g$ is
called {\it type-rotating\/}. Such automorphisms form a subgroup of
index at most~2 in the group of all automorphisms of~$\Delta$.

If $K$ is a field with discrete valuation, then there is a thick
building $\Delta_K$ of type $\widetilde {A_n}$ associated with $K$
[\Ro, Section 9.2], [\B, Section V.8], and the group $PGL(n+1, K)$
acts transitively and in a type-rotating way on $\Delta_K$.

A group is said to be an {\it $\widetilde A_n$-group} if it acts
simply transitively on the vertices of a thick building of type
$\widetilde A_n$ in a type-rotating way.

We now describe $\widetilde{A_n}$-groups. Let $\Pi$ be a projective
geometry of dimension~$n\ge2$ (see [\D, p.~24] or [\T, p.~105], for
example).  For $i=1,\ldots,n$, let $\Pi_i=\{x\in\Pi:\dim(x)=i\}$. To
avoid unnecessary abstraction, the reader may assume that $\Pi$ is the
set $\Pi(\V)$ (partially ordered by inclusion) of nontrivial proper
subspaces of an $n+1$-dimensional vector space~$\V$ over a field~$k$,
and that $\dim(x)$ refers to the dimension of the subspace $x$
of~$\V$. For when $n\ge3$, or when $n=2$ and $\Pi$ is desarguesian,
$\Pi$ must be isomorphic to~$\Pi(\V)$ for some~$\V$ ([\D, pp.~27--28]
or [\T, p.~203]). Let $\lambda:\Pi\to\Pi$ be an involution such that
$\lambda(\Pi_i)=\Pi_{n+1-i}$ for $i=1,\ldots,n$, and let $\cT$ is an
{\it $\widetilde{A_n}$-triangle presentation\/} compatible
with~$\lambda$. This means that $\cT$ is a set of triples $(u,v,w)$,
where $u,v,w\in\Pi$, such that

\smallskip

\item{(A)} given $u,v\in\Pi$, then $(u,v,w)\in\cT$ for some $w\in\Pi$
if and only if $\lambda(u)$ and~$v$ are distinct and incident;
\smallskip

\item{(B)} if $(u,v,w)\in\cT$, then $(v,w,u)\in\cT$;
\smallskip

\item{(C)} if $(u,v,w_1)\in\cT$ and $(u,v,w_2)\in\cT$, then $w_1=w_2$;
\smallskip

\item{(D)} if $(u,v,w)\in\cT$, then
$(\lambda(w),\lambda(v),\lambda(u))\in\cT$; \smallskip

\item{(E)} if $(u,v,w)\in\cT$, then $\dim(u)+\dim(v)+\dim(w)= n+1$
or~$2(n+1)$; \smallskip
   
\item{(F)} if $(x,y,u)\in\cT'$ and $(x',y',\lambda(u))\in\cT'$, then
for some $w\in\Pi$ we have $(y',x,w)\in\cT'$ and
$(y,x',\lambda(w))\in\cT'$.  \smallskip
 
\noindent Here $\cT'$ denotes the ``half" of~$\cT$ consisting of the
triples $(u,v,w)\in\cT$ for which $\dim(u)+\dim(v)+\dim(w)= n+1$.
Given $u,v\in\Pi$, then $(u,v,w)\in\cT'$ for some $w\in\Pi$ if and
only if $\lambda(u)\supsetneqq v$. We also write $\cT''$ for
$\cT\setminus\cT'$.

We form the associated group $\Gamma_\cT$ with a generating set
indexed by~$\Pi$:
$$ 
\eqalign{
\Gamma_\cT=\langle\{a_v\}_{v\in\Pi}\mid
&\ (1)\ a_{\lambda(v)}=a_v^{-1}{\rm\ for\ all\ }v\in\Pi,\cr
&\ (2)\ a_ua_va_w=1{\rm\ for\ all\ }(u,v,w)\in\cT\rangle.\cr}
$$ 
It was shown in~[\C] that the Cayley graph of~$\Gamma_\cT$ with
respect to the generators $a_v$, $v\in\Pi$, is the 1-skeleton of a
thick building~$\Delta_\cT$ of type~$\widetilde{A_n}$. Clearly,
$\Gamma_\cT$ acts, by left multiplication, simply transitively on the
set of vertices of~$\Delta_\cT$. Conversely, if $\Gamma$ is a group of
type-rotating automorphisms of a thick building~$\Delta$ of
type~$\widetilde{A_n}$, and acts simply transitively on the vertices
of~$\Delta$, then $\Gamma\cong\Gamma_\cT$ and $\Delta\cong\Delta_\cT$
for some $\widetilde{A_n}$-triangle presentation~$\cT$. This
generalized earlier work on the case $n=2$ [\CMSZa].

In this paper, $\Pi$ is assumed finite. The number of $x\in\Pi_1$
incident with any given $y\in\Pi_2$ is denoted $q+1$, and is
independent of~$y$. Here $q$ is called the {\it order\/} of~$\Pi$, and
when $\Pi=\Pi(\V)$, $q$ is the number of elements in the field~$k$. It
remains an open problem whether $\widetilde{A_n}$-triangle
presentations exist for every~$n$. They have been found when
$\Pi=\Pi(\V)$ for $n=2$ [\CMSZa] and for $n=3,4$ [\C] for any prime
power~$q$, and for $n=5$ and $q=2$ [Cartwright, unpublished].

Let $L$ denote the set of all strings $u_1u_2\cdots u_\ell$ over~$\Pi$
such that $\lambda(u_i)+u_{i+1}=\V$ for $i=1,\ldots,\ell-1$. (The
notation assumes that $\Pi=\Pi(\V)$, but in general,
``$\lambda(u_i)+u_{i+1}=\V$" is interpreted as ``there is no $x\in\Pi$
such that $\lambda(u_i)\subset x$ and $u_{i+1}\subset x$, where we
write $y\subset x$ if $x$ and~$y$ are incident, and
$\dim(y)\le\dim(x)$.) Theorem~2.2 in~[\C] states that $u=u_1\cdots
u_\ell\mapsto {\bar u}=a_{u_1}\cdots a_{u_\ell}$ is a bijection
$L\to\Gamma_\cT$. Moreover, in the notation of~[\Ep], each string
$u_1\cdots u_\ell\in L$ is {\it geodesic\/}, and so the number~$\ell$
is the {\it word length\/} $|g|$ of~$g={\bar u}$. If $g\in\Gamma_\cT$,
$u\in L$ and $g={\bar u}$, then $u$ is called the {\it normal form\/}
of~$g$. We shall also refer to strings in~$L$ as being {\it in normal
form\/}.  When $g$ is the identity element~1, its normal form is the
empty word, and its word length is~0, by definition. We write ${\rm
d}(g,g')$ for $|g^{-1}g'|$, the distance from~$g$ to~$g'$ in the {\it
word metric\/}.

\medskip\noindent{\bf\S5. Finitely generated $\widetilde {A_n}$-groups 
are automatic.} \medskip 

\noindent Throughout this section, let $\Gamma=\Gamma_\cT$ be a
finitely generated $\widetilde{A_n}$-group.

We start by showing that  $L$ is a regular language. For we can define
a finite state automaton~$M$ accepting~$L$ as follows: the set~$S$ of
states of~$M$ is $\Pi$ together with an initial state $s_0$ (not
in~$\Pi$) and a single failure state $s_1$ (not in~$\Pi$, and distinct
from~$s_0$); the alphabet of~$M$ is~$\Pi$; the transition function
$\mu$ of~$M$ is given by $\mu(s_0,x)=x$ and $\mu(s_1,x)=s_1$ for
$x\in\Pi$, while for $x,y\in\Pi$, we set $\mu(x,y)=y$ if
$\lambda(x)+y=\V$ and $\mu(x,y)=s_1$ otherwise; the set~$Y$ of accept
states of~$M$ is $\Pi\cup\{s_0\}$. Clearly $L=L(M)$.

Note that $\Pi$ is a set of semigroup generators for $\Gamma_\cT$
which is closed under inversion (as $a_{\lambda(x)}=a_x^{-1}$).
Moreover, $L$ has the {\it uniqueness property\/} (i.e.,
$u\mapsto{\bar u}$ is a bijection $L\to\Gamma_\cT$), is obviously
prefix-closed, and is symmetric (i.e., $u_1\cdots u_\ell\in L$ implies
that $\lambda(u_\ell)\cdots\lambda(u_1)\in L$).  This last property
of~$L$ and our Theorem~2 below imply that $\Gamma_\cT$ is {\it
biautomatic\/}, and in fact {\it symmetric automatic} or {\it fully
automatic\/}.

 If $u=u_1\cdots u_\ell\in L$, let $u(0)$ be
the identity element~1 in~$\Gamma$, and for $1\le t\le\ell$, let
$u(t)$ denote the element $a_{u_1}\cdots a_{u_t}$ of~$\Gamma$. For
$t>\ell$, let $u(t)=a_{u_1}\cdots a_{u_\ell}={\bar u}$. Let $x\in\Pi$,
and let $v_1\cdots v_k$ be the normal form
of~$\overline{ux}=a_{u_1}\cdots a_{u_\ell}a_x$. By [\Ep,
Theorem~2.3.5], to show that $\Gamma_\cT$ is automatic, it is enough
to show that for some $k$ (independent of~$u$ and~$x$),
$d(u(t),v(t))\le k$ holds for each integer $t\ge0$. This property is
called the {\it $k$-fellow traveller property} [\BGSS]. In fact, we
show this property holds for $k=1$.

Now let $u=u_1\cdots u_\ell\in L$, and let $x\in\Pi$. The main work
below is to describe the normal form $v=v_1\cdots v_k$
of~$a_{u_1}\cdots a_{u_\ell} a_x$. The description of this normal form
is complicated by the fact that, given $u,x\in\Pi$, there are the
following~5 mutually exclusive cases to consider: $\lambda(u)=x$,
$\lambda(u)+x=\V$, $\lambda(u)\supsetneqq x$, $\lambda(u)\subsetneqq
x$, and, finally, $\lambda(u)$ and~$x$ are distinct and nonincident,
with $\lambda(u)+x\ne\V$.

We start with a lemma which concerns the last of these cases. In this
lemma and below, we consistently identify a string~$u$ over~$\Pi$ with
its image $\bar u$ in~$\Gamma_\cT$, sometimes writing ``$u$ in
$\Gamma_\cT$" for emphasis, when thinking of~$u$ as~$\bar u$.

\medskip\noindent{\bf Lemma 1.}\enspace{\sl Suppose that $u,x\in\Pi$, and that
$\lambda(u)$ and~$x$ are distinct and nonincident, and that 
$\lambda(u)+x=\lambda(s)$ ($\ne\V$). Thus we can write 
$(s,\lambda(u),x')\in\cT'$ and $(s,x,\lambda(v))\in\cT'$
for some $x',v\in\Pi$, and $ux=x'v$ in~$\Gamma$. Then $x'v$ is in normal form.
That is, $\lambda(x')+v=\V$. Conversely,
if two triples $(s,\lambda(u),x'),(s,x,\lambda(v))\in\cT'$ are given, with 
$\lambda(x')+v=\V$, then $\lambda(u)$ and~$x$ must be distinct and nonincident,
with $\lambda(u)+x=\lambda(s)$.}\par
\medskip
We illustrate this lemma with a diagram in the Cayley graph of~$\Gamma_\cT$: 
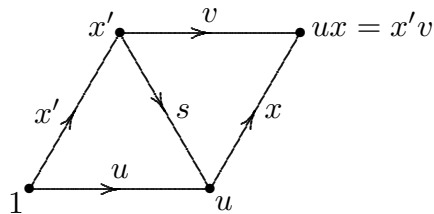
$$
\beginpicture
\setcoordinatesystem units <1mm,1mm>
\setplotarea x from -20 to 50, y from -5 to 20
\arrow <6pt> [.15,.6] from 0 0 to 10 0
\arrow <6pt> [.15,.6] from 0 0 to 5 8.66
\arrow <6pt> [.15,.6] from 10 17.32 to 20 17.32
\arrow <6pt> [.15,.6] from 10 17.32 to 15 8.66
\arrow <6pt> [.15,.6] from 20 0 to 25 8.66
\setlinear 
\plot 5 8.66 10 17.32 /
\plot 10 0 20 0 15 8.66 /
\plot 25 8.66 30 17.32 /
\plot 20 17.32 30 17.32 / 
\put {$x'$} at 2 8.66
\put {$u$} at 10 2
\put {$s$} at 17 8.66
\put {$v$} at 20 19.32
\put {$x$} at 27 8.66  
\put {$\scriptstyle\bullet$} at 0 0
\put {$\scriptstyle\bullet$} at 20 0
\put {$\scriptstyle\bullet$} at 10 17.32
\put {$\scriptstyle\bullet$} at 30 17.32 
\put {1} at -1.5 -1.5
\put {$u$} at 21.5 -1.5
\put {$x'$} at 8 18.32
\put {$ux=x'v$} at 38 18.32
\endpicture 
$$
\medskip\noindent{\bf Proof.} 
Firstly, $\lambda(x')\ne v$. Otherwise, 
$x'=\lambda(v)$, and so Axioms~(B) and~(C) imply that $\lambda(u)=x$,
contrary to hypothesis. We shall henceforth use Axioms (A)--(E) in the
definition of an $\widetilde{A_n}$-triangle presentation without comment,
but refer to Axiom~(F) when it is used.
Next, suppose that $\lambda(x')\ne v$, but that $\lambda(x')$ and~$v$ are 
incident. Thus $(x',v,w)\in\cT$ for some $w\in\Pi$.
If $(x',v,w)\in\cT'$, then $(w,x',v)\in\cT'$, 
$(s,x,\lambda(v))\in\cT'$ and Axiom~(F) imply that $(x,w,y)\in\cT'$ and
$(x',s,\lambda(y))\in\cT'$ for some $y\in\Pi$. Thus $y=u$. But then $(u,x,w)\in\cT$, so that $\lambda(u)$ and~$x$ are
incident, contrary to hypothesis.   
If $(x',v,w)\in\cT''$, then $(s,\lambda(u),x')\in\cT'$,
$(\lambda(w),\lambda(v),\lambda(x'))\in\cT'$ and Axiom~(F) imply that
$(\lambda(v),s,y)\in\cT'$ and $(\lambda(u),\lambda(w),\lambda(y))\in\cT'$
for some $y\in\Pi$. Then $y=x$, and so
$(x,w,u)\in\cT''$, so that again $\lambda(u)$ and~$x$ are
incident, contrary to hypothesis.   

Suppose that $\lambda(x')\ne v$, that $\lambda(x')$ and~$v$ are not
incident, and that $\lambda(x')+v=\lambda(s')\ne\V$. So 
$(s',\lambda(x'),x'')\in\cT'$ and $(s',v,\lambda(v'))\in\cT'$ for some 
$x'',v'\in\Pi$. Now $(s,\lambda(u),x')\in\cT'$, $(x'',s',\lambda(x'))\in\cT'$ and
Axiom~(F) show that $(s',s,z)\in\cT'$ and 
$(\lambda(u),x'',\lambda(z))\in\cT'$ for some $z\in\Pi$. Similarly,
$(\lambda(v'),s',v)\in\cT'$ and $(s,x,\lambda(v))\in\cT'$ imply that
$(x,\lambda(v'),z')\in\cT'$ and $(s',s,\lambda(z'))\in\cT'$ for some 
$z'\in\Pi$. Thus $\lambda(z')=z$. Hence
$(\lambda(u),x'',\lambda(z))\in\cT'$ and 
$(x,\lambda(v'),\lambda(z))\in\cT'$, so that $z\supset\lambda(u),x$.
Hence $z\supset\lambda(u)+x=\lambda(s)$. But $(s',s,z)\in\cT'$
implies that $z\subsetneqq\lambda(s)$. This contradiction completes
the proof of the first part of the lemma.  

Consider the converse part. Firstly, $\lambda(u)$ and~$x$ must
be distinct, for otherwise $x'=\lambda(v)$ must hold, which is impossible,
because $\lambda(x')+v=\V$. Next, $\lambda(u)\supsetneqq x$ cannot hold. For
otherwise $(u,x,y)\in\cT'$ for some $y\in\Pi$. This, 
$(\lambda(u),x',s)\in\cT'$ and Axiom~F then imply that $(s,x,z)\in\cT'$ and
$(y,x',\lambda(z))\in\cT'$ hold for some $z\in\Pi$. Thus $z=\lambda(v)$,
so that $(y,x',v)\in\cT'$, which implies that $\lambda(x')\supset v$,
again a contradiction. Similarly, $\lambda(u)\subsetneqq x$ leads to a 
contradiction. Finally, the hypotheses imply that 
$\lambda(u)+x\subset\lambda(s)$. If $\lambda(u)+x\ne\lambda(s)$, write
$\lambda(u)+x=\lambda(s')$, and then we can find $x'',v'\in\Pi$ such that
$(s',\lambda(u),x''),(s',x,\lambda(v'))\in\cT'$. But then $x'v=ux=x''v'$,
and both $x'v$ and $x''v'$ are in normal form (by the hypotheses, and by
the first part of the lemma, respectively). Uniqueness of normal forms now
shows that $x''=x'$, so that $s'=s$, a contradiction. This completes the proof.

\medskip
\medskip\noindent{\bf Lemma 2.}\enspace{\sl Let $u_1,u_2,x\in\Pi$, with
$u_1u_2$ in normal form. Then to get the normal form of $u_1u_2x$, there are
the following 5~possibilities:
\smallskip
\item{(1)} If $\lambda(u_2)+x=\V$, then $u_1u_2x$ is in normal form.
\smallskip
\item{(2)} If $\lambda(u_2)=x$, then the normal form of $u_1u_2x$ is~$u_1$.
\smallskip
\item{(3)} If $\lambda(u_2)\supsetneqq x$, then $(u_2,x,\lambda(w))\in\cT'$ for
some $w\in\Pi$, and the normal form of $u_1u_2x$ is~$u_1w$.
\smallskip
\item{(4)} If $\lambda(u_2)\subsetneqq x$, then $(u_2,x,\lambda(w))\in\cT''$ for
some $w\in\Pi$. Thus $u_1u_2x=u_1w$ in~$\Gamma$. There are
now the following possibilities: 
\smallskip
\itemitem{(a)} Either $\lambda(u_1)+w=\V$, in which case the normal form of
$u_1u_2x$ is~$u_1w$, or 
\smallskip
\itemitem{(b)} $\lambda(u_1)\supsetneqq w$, in which case 
$(u_1,w,\lambda(w'))\in\cT'$ for
some $w'\in\Pi$, and the normal form of $u_1u_2x$ is $w'$, or  
\smallskip
\itemitem{(c)} $\lambda(u_1)$ and $w$ are distinct and nonincident, with  
$\lambda(u_1)+w\ne\V$. Then writing $\lambda(u_1)+w=\lambda(s')$, there are
unique $w',v_1\in\Pi$ such that $(s',\lambda(u_1),w')\in\cT'$ and 
$(s',w,\lambda(v_1))\in\cT'$, and the normal form of $u_1u_2x$ is $w'v_1$.
\smallskip
\item{(5)} If $\lambda(u_2)$ and $x$ are distinct and nonincident, with  
$\lambda(u_2)+x\ne\V$, then writing $\lambda(u_2)+x=\lambda(s)$, there are
unique $x',v_2\in\Pi$ such that $(s,\lambda(u_2),x')\in\cT'$ and 
$(s,x,\lambda(v_2))\in\cT'$. Thus $u_1u_2x=u_1x'v_2$ in~$\Gamma$. There are
now the following possibilities: 
\smallskip
\itemitem{(a)} Either $\lambda(u_1)+x'=\V$, in which case the normal form of
$u_1u_2x$ is~$u_1x'v_2$, or 
\smallskip
\itemitem{(b)} $\lambda(u_1)\supsetneqq x'$, in which case 
$(u_1,x',\lambda(w'))\in\cT'$ for
some $w'\in\Pi$, and the normal form of $u_1u_2x$ is $w'v_2$, or
\smallskip
\itemitem{(c)} $\lambda(u_1)$ and $x'$ are distinct and nonincident, with  
$\lambda(u_1)+x'\ne\V$. Then writing $\lambda(u_1)+x'=\lambda(s')$, there are
unique $x'',v_1\in\Pi$ such that $(s',\lambda(u_1),x'')\in\cT'$ and 
$(s',x',\lambda(v_1))\in\cT'$, and the normal form of $u_1u_2x$ is $x''v_1v_2$.
\par} 
We illustrate the more complicated cases:
$$
\beginpicture 
\setcoordinatesystem units <1mm,1mm> point at 45 8.66
\setplotarea x from -10 to 57, y from -5 to 20
\arrow <6pt> [.15,.6] from 0 0 to 10 0
\arrow <6pt> [.15,.6] from 0 0 to 5 8.66
\arrow <6pt> [.15,.6] from 10 17.32 to 20 17.32
\arrow <6pt> [.15,.6] from 10 17.32 to 15 8.66
\arrow <6pt> [.15,.6] from 20 0 to 25 8.66
\arrow <6pt> [.15,.6] from 20 0 to 30 0
\arrow <6pt> [.15,.6] from 40 0 to 35 8.66
\setlinear 
\plot 5 8.66 10 17.32 /
\plot 10 0 20 0 15 8.66 /
\plot 25 8.66 30 17.32 20 17.32 /
\plot 30 0 40 0 / 
\plot 35 8.66 30 17.32 / 
\put {$u_1$} at 10 2
\put {$s'$} at 17.5 8.66
\put {$v_1$} at 20 15.32
\put {$u_2$} at 30 2 
\put {$x$} at 37 8.66  
\put {$w$} at 27 8.66  
\put {$w'$} at 2 8.66  
\put {$\scriptstyle\bullet$} at 0 0
\put {$\scriptstyle\bullet$} at 20 0
\put {$\scriptstyle\bullet$} at 40 0
\put {$\scriptstyle\bullet$} at 10 17.32
\put {$\scriptstyle\bullet$} at 30 17.32 
\put {1} at -1.5 -1.5
\put {$u_1$} at 21.5 -2
\put {$u_1u_2$} at 41.5 -2 
\put {$w'$} at 8 18.32
\put {$u_1u_2x=w'v_1$} [lb] at 31.5 16.32
\put {Case 4(c)} at 20 -8
\endpicture 
\beginpicture 
\setcoordinatesystem units <1mm,1mm> point at 10 8.66
\setplotarea x from -10 to 57, y from -5 to 20
\arrow <6pt> [.15,.6] from -10 17.32 to -5 8.66
\arrow <6pt> [.15,.6] from -10 17.32 to 0 17.32
\arrow <6pt> [.15,.6] from 0 0 to 5 8.66
\arrow <6pt> [.15,.6] from 0 0 to 10 0
\arrow <6pt> [.15,.6] from 10 17.32 to 20 17.32
\arrow <6pt> [.15,.6] from 10 17.32 to 15 8.66
\arrow <6pt> [.15,.6] from 20 0 to 25 8.66
\setlinear 
\plot 5 8.66 10 17.32 0 17.32 /
\plot -5 8.66 0 0 /
\plot 10 0 20 0 15 8.66 / 
\plot 25 8.66 30 17.32 20 17.32 / 
\put {$u_1$}  at -7.5 7.66
\put {$s$} at 17 8.66
\put {$x'$} at 7 8.66  
\put {$u_2$} at 10 -2 
\put {$x$} at 27 8.66   
\put {$w'$} at 0 15.32
\put {$v_2$} at 21 15.32
\put {$\scriptstyle\bullet$} at 0 0
\put {$\scriptstyle\bullet$} at 20 0
\put {$\scriptstyle\bullet$} at -10 17.32
\put {$\scriptstyle\bullet$} at 10 17.32 
\put {$\scriptstyle\bullet$} at 30 17.32 
\put {1} at -12 17.32
\put {$u_1$} at -1 -2
\put {$u_1x'=w'$} at 10 19.32
\put {$u_1u_2$} at 21 -2 
\put {$u_1u_2x=w'v_2$} [lb] at 31.5 16.32
\put {Case 5(b)} at 10 -8
\endpicture  
$$
$$
\beginpicture 
\setcoordinatesystem units <1mm,1mm>
\setplotarea x from -10 to 80, y from -5 to 20
\arrow <6pt> [.15,.6] from 0 0 to 10 0
\arrow <6pt> [.15,.6] from 0 0 to 5 8.66
\arrow <6pt> [.15,.6] from 10 17.32 to 20 17.32
\arrow <6pt> [.15,.6] from 10 17.32 to 15 8.66
\arrow <6pt> [.15,.6] from 20 0 to 25 8.66
\arrow <6pt> [.15,.6] from 20 0 to 30 0
\arrow <6pt> [.15,.6] from 30 17.32 to 42 17.32
\arrow <6pt> [.15,.6] from 30 17.32 to 35 8.66
\arrow <6pt> [.15,.6] from 40 0 to 45 8.66
\setlinear 
\plot 5 8.66 10 17.32 /
\plot 10 0 20 0 15 8.66 /
\plot 25 8.66 30 17.32 /
\plot 20 17.32 30 17.32 / 
\plot 30 0 40 0 35 8.66 / 
\plot 45 8.66 50 17.32 40 17.32 / 
\put {$x''$} at 1 8.66
\put {$u_1$} at 10 2
\put {$s'$} at 17.5 8.66
\put {$v_1$} at 20 15.32
\put {$x'$} at 27 8.66  
\put {$u_2$} at 30 2 
\put {$s$} at 37 8.66
\put {$x$} at 47 8.66  
\put {$v_2$} at 41 15.32
\put {$\scriptstyle\bullet$} at 0 0
\put {$\scriptstyle\bullet$} at 20 0
\put {$\scriptstyle\bullet$} at 40 0
\put {$\scriptstyle\bullet$} at 10 17.32
\put {$\scriptstyle\bullet$} at 30 17.32 
\put {$\scriptstyle\bullet$} at 50 17.32 
\put {1} at -1.5 -1.5
\put {$u_1$} at 21.5 -2
\put {$x''$} at 8 18.32
\put {$u_1x'=x''v_1$} at 29 19.32
\put {$u_1u_2$} at 41.5 -2 
\put {$u_1u_2x=x''v_1v_2$} [lb] at 51.5 16.32
\put {Case 5(c)} at 30 -8
\endpicture
$$

\medskip\noindent{\bf Proof.} The assertions in~(1) and~(2) are obvious. 

Let us consider the situation in~(3).  Then $(u_2,x,\lambda(w))\in\cT'$
implies that $w\supset u_2$. Thus $\lambda(u_1)+w\supset\lambda(u_1)+u_2=\V$,
and so $u_1w$ is in normal form.

Let us consider the situation in~(4). In this case, again $w\ne\lambda(u_1)$. 
For otherwise, $(u_2,x,u_1)=(u_2,x,\lambda(w))\in\cT''$, contradicting the hypothesis that $u_1u_2$ is in 
normal form. Also, $\lambda(u_1)\subsetneqq w$ cannot hold. For otherwise,  
$(\lambda(w),\lambda(u_1),w')\in\cT'$ for some $w'\in\Pi$. Now
$(\lambda(x),\lambda(u_2),w)\in\cT'$, 
$(\lambda(u_1), w', \lambda(w))\in\cT'$ and Axiom~(F) imply that
$(w',\lambda(x),y)$ and 
$(\lambda(u_2),\lambda(u_1), \lambda(y))$ are in~$\cT'$ for
some $y\in\Pi$. This again contradicts the hypothesis that $u_1u_2$ is in 
normal form. We are left with the three possibilities 4(a), 4(b) and~4(c). The
assertions in~4(a) and~4(b) are obvious, while in~4(c), the fact that 
the word $w'v_1$ is in normal form is immediate from Lemma~1. 

Let us consider the situation in~(5). In this case, $\lambda(u_1)$
cannot equal~$x'$. For otherwise $(s,\lambda(u_2),\lambda(u_1))=
(s,\lambda(u_2),x')\in\cT'$, which contradicts the hypothesis 
that $u_1u_2$ is in normal form. Also, $\lambda(u_1)\subsetneqq x'$
cannot happen. For otherwise, 
$(\lambda(x'),\lambda(u_1), w')\in\cT'$ for some $w'\in\Pi$. Now
$(s,\lambda(u_2),x')\in\cT'$, 
$(\lambda(u_1), w',\lambda(x'))\in\cT'$ and Axiom~(F) imply that
$(w',s,y)\in\cT'$ and 
$(\lambda(u_2),\lambda(u_1),\lambda(y))\in\cT'$ for some $y\in\Pi$.
This last fact again contradicts the hypothesis that $u_1u_2$ is in
normal form. We are left with the three possibilities 5(a), 5(b) 
and~5(c). The assertion in~5(a) is obvious, by Lemma~1.

Consider the situation in~5(b). First notice that $(x',\lambda(w'),u_1),
(x',s,\lambda(u_2))\in\cT'$ and the converse part of Lemma~1 imply that
$\lambda(w')+s=\lambda(x')$. As $(s,x,\lambda(v_2))\in\cT'$, we have 
$v_2\supset s$, and so $\lambda(w')+v_2=\lambda(w')+s+v_2=\lambda(x')+v_2=\V$,
the last equation holding by Lemma~1. So $w'v_2$ is in normal form.

Finally, consider the situation in~5(c). First observe that 
$\lambda(v_1)+s=\lambda(x')$. This
follows from the converse part of Lemma~1, because 
$(x',\lambda(v_1),s'),(x',s,\lambda(u_2))\in\cT'$, and because, using
$(s',\lambda(u_1),x'')\in\cT'$, we have 
$\lambda(s')+u_2\supset\lambda(u_1)+u_2=\V$. Now $(s,x,\lambda(v_2))\in\cT'$,
so that $v_2\supset s$. Hence $\lambda(v_1)+v_2=\lambda(v_1)+s+v_2
=\lambda(x')+v_2=\V$, the last equation holding by Lemma~1. Lemma~1 also
shows that $\lambda(x'')+v_1=\V$, and so $x''v_1v_2$ is in normal form.

\medskip\noindent{\bf Remark.}\enspace In later work, we shall need a converse
to Lemma~2. Let us write $x+'y=z$ for $x,y,z\in\Pi$ if $x,y$ are distinct and
nonincident, with $x+y=z$. Then by the converse part of Lemma~1, we have
$\lambda(w')+'\lambda(x)=\lambda(w)$ in part~4(b),
$\lambda(v_1)+'\lambda(x)=\lambda(w)$ in part~4(c),
$\lambda(w')+'s=\lambda(x')$ in part~5(b) and
$\lambda(v_1)+'s=\lambda(x')$ in part~5(c). Provided these conditions are 
added, the converses of parts~(4) and~(5) hold. For example, for 5(b), if 
$u_1,u_2,x,s,x',v_2,w'\in\Pi$ and triples 
$(s,\lambda(u_2),x'),(s,x,\lambda(v_2))\in\cT'$ are given with 
$\lambda(w')+'s=\lambda(x')$ and $\lambda(w')+v_2=\V$, then $u_1u_2$ is in
normal form and $\lambda(u_2)+'x=\lambda(s)$. The first of these is immediate 
from Lemma~1, since $(x',\lambda(w'),u_1),(x',s,\lambda(u_2))\in\cT'$. To see
that $\lambda(u_2)+'x=\lambda(s)$, notice that 
$\lambda(x')+v_2\supset\lambda(w')+v_2=\V$, and so the converse part of Lemma~1
is applicable.

\medskip\noindent{\bf The normal form $v_1\cdots v_k$ 
of $a_{u_1}\cdots a_{u_\ell}a_x$.}
\medskip
We can now describe how to obtain the normal form $v_1\cdots v_k$ 
of $\overline{ux}=a_{u_1}\cdots a_{u_\ell}a_x$, given $u_1\cdots u_\ell\in L$ 
and $x\in\Pi$. When $\lambda(u_\ell)+x=\V$, this is obviously 
$u_1\cdots u_\ell x$. If $\lambda(u_\ell)=x$, then the normal 
form of $\overline{ux}$ is clearly $u_1\cdots u_{\ell-1}$. When 
$\lambda(u_\ell)\supsetneqq x$, then $(u_\ell,x,\lambda(w))\in\cT'$ for
some $w\in\Pi$, and 
the normal form of $\overline{ux}$ is $u_1\cdots u_{\ell-1}w$, by
part~(3) of Lemma~2.

When $\lambda(u_\ell)$ and $x$ are distinct and nonincident, with
$\lambda(u_\ell)+x\ne\V$, then writing $\lambda(u_\ell)+x=\lambda(s_\ell)$,
we have $(s_\ell,\lambda(u_\ell),x_1)\in\cT'$ and 
$(s_\ell,x,\lambda(v_\ell))\in\cT'$ for some $x_1,v_\ell\in\Pi$. Then
$u_1\cdots u_\ell x=u_1\cdots u_{\ell-1}x_1v_\ell$ in~$\Gamma$. This
situation may be repeated several times, with $\lambda(u_{\ell-1})$ 
and~$x_1$ distinct and nonincident, with $\lambda(u_{\ell-1})+x_1\ne\V$,
and so on. Suppose this situation is repeated exactly $i$~times. Then------ end of 1994.5 -- ascii -- complete ------

we find $s_{\ell-\nu}, x_{\nu+1}, v_{\ell-\nu}\in\Pi$ for $\nu=0,\ldots,i-1$,
such that $\lambda(u_{\ell-\nu})+x_\nu=\lambda(s_{\ell-\nu})$,
$(s_{\ell-\nu},\lambda(u_{\ell-\nu}),x_{\nu+1})\in\cT'$, and
$(s_{\ell-\nu}, x_\nu,\lambda(v_{\ell-\nu}))\in\cT'$ for $\nu=0,\ldots,i-1$
(writing $x_0=x$). Then in~$\Gamma$,
$$
u_1\cdots u_\ell x=u_1\cdots u_{\ell-i}x_iv_{\ell-i+1}\cdots v_\ell\,.\eqno(4.1)
$$ 
Then by part~(5) of Lemma~2, either the word on the right in~(4.1)
is in normal form, or $\lambda(u_{\ell-i})\supsetneqq x_i$, in which case
$(u_{\ell-i},x_i,\lambda(w))\in\cT'$ for some $w\in\Pi$, and the normal form
of~$\overline{ux}$ is
$$
u_1\cdots u_{\ell-i-1}wv_{\ell-i+1}\cdots v_\ell\,.\eqno(4.2)
$$ 

Finally, if $\lambda(u_\ell)\subsetneqq x$, then $(u_\ell,x,\lambda(w))\in\cT''$
for some $w\in\Pi$. Thus $u_1\cdots u_\ell x=u_1\cdots u_{\ell-1}w$ in~$\Gamma$.
A sequence of exactly $i\ge0$ steps such as led to~(4.1) may now occur,
so that we find elements $s_{\ell-\nu}, w_\nu, v_{\ell-\nu}$ for 
$\nu=1,\ldots,i$, such that 
$\lambda(u_{\ell-\nu})+w_{\nu-1}=\lambda(s_{\ell-\nu})$,
$(s_{\ell-\nu},\lambda(u_{\ell-\nu}),w_\nu)\in\cT'$ and
$(s_{\ell-\nu}, w_{\nu-1},\lambda(v_{\ell-\nu}))\in\cT'$ for $\nu=1,\ldots,i$
(writing $w_0=w$). Then in~$\Gamma$,
$$
u_1\cdots u_\ell x=u_1\cdots u_{\ell-i-1}w_iv_{\ell-i}\cdots v_{\ell-1}\,.\eqno(4.3)
$$ 
By parts~(4) and~(5) of Lemma~2, either the word on the right in~(4.3)
is in normal form, or $\lambda(u_{\ell-i-1})\supsetneqq w_i$, in which case
$(u_{\ell-i-1},w_i,\lambda(w'))\in\cT'$ for some $w'\in\Pi$, and the normal form
of~$\overline{ux}$ is
$$
u_1\cdots u_{\ell-i-2}w'v_{\ell-i}\cdots v_{\ell-1}\,.\eqno(4.4)
$$ 

We are now ready to prove the fellow traveller property.
\medskip\noindent{\bf Theorem 2.}\enspace{\sl Let $\Gamma_\cT$ be a
finitely generated $\widetilde{A_n}$-group. Let $u=u_1\cdots u_\ell\in L$, 
let $x\in\Pi$, and let $v_1\cdots v_k\in L$ be the normal form of 
$\overline{ux}=a_{u_1}\cdots a_{u_\ell}a_x$. Then for each integer $t\ge0$, 
we have ${\rm d}(u(t),v(t))\le1$. That is, either $a_{u_1}\cdots
a_{u_t}=a_{v_1}\cdots a_{v_t}$, or $a_{u_1}\cdots
a_{u_t}=a_{v_1}\cdots a_{v_t}a_{x_t}$ for some $x_t\in\Pi$. Thus
$\Gamma_\cT$ is an automatic group.}\par

\medskip\noindent{\bf Proof.} Let us discuss the most complicated case
in detail, leaving the other cases to the reader. Suppose that
$\lambda(u_\ell)\subsetneqq x$, and that the normal form~$v$
of~$\overline{ux}$ is~(4.4). For $0\le t\le\ell-i-2$, we have
$a_{u_1}\cdots a_{u_t}=a_{v_1}\cdots a_{v_t}$. Now
$$
\eqalign{v(\ell-i-1)=a_{v_1}\cdots a_{v_{\ell-i-1}}
&=a_{u_1}\cdots a_{u_{\ell-i-2}}a_{w'}\cr
&=a_{u_1}\cdots a_{u_{\ell-i-2}}a_{u_{\ell-i-1}}a_{w_i}\cr
&=u(\ell-i-1)a_{w_i}\cr}
$$
and $w_i\in\Pi$. For $1\le\nu\le i$, we have
$$
\eqalign{v(\ell-\nu)
&=a_{u_1}\cdots a_{u_{\ell-i-2}}a_{w'}a_{v_{\ell-i}}\cdots a_{v_{\ell-\nu}}\cr
&=a_{u_1}\cdots a_{u_{\ell-i-2}}a_{u_{\ell-i-1}}a_{w_i}a_{v_{\ell-i}}
\cdots a_{v_{\ell-\nu}}\cr
&=a_{u_1}\cdots a_{u_{\ell-i-2}}a_{u_{\ell-i-1}}
         a_{u_{\ell-i}}a_{u_{\ell-i-1}}\cdots a_{u_{\ell-\nu}}a_{w_{\nu-1}}\cr 
&=u(\ell-\nu)a_{w_{\nu-1}}\cr}
$$
and $w_{\nu-1}\in\Pi$. For $t\ge\ell$, $v(t)=\overline{ux}=
u(t)a_x$. This completes the proof.  \bigskip 

\centerline{\bf References.}  \medskip

\item{[\ABC]} Alonso, J., Brady, T., Cooper, D., Ferlini, V., Lustig,
M., Mihalik, M., Shapiro, M., and Short, H., Short, H., ed., Notes on
word hyperbolic groups in ``Group Theory From a Geometric Viewpoint'',
Ghys, E., Haefliger, A., and Verjovsky A., eds., World Scientific (1990).
\smallskip

\item{[\BGSS]} Baumslag, G., Gersten, S.M., Shapiro, M. and Short, H.,
Automatic groups and amalgams, J. Pure and Appl. Alg. {\bf76} (1991), 229--316.
\smallskip

\item{[\B]} Brown, K., ``Buildings'', Springer Verlag (1989).
\smallskip

\item{[\Can]} Cannon, J.W., The combinatorial structure of cocompact
discrete hyperbolic groups, {\it Geom. Ded.\/}, {\bf 16} (1984), 123--148.

\item{[\C]} Cartwright, D.I., Groups acting simply transitively on the 
vertices of a building of type~${\tilde A}_n$, to appear, Proceedings
of the 1993 Como conference ``Groups of Lie type and their
geometries'', W.M. Kantor, Editor.
\smallskip

\item{[\CMSZa]} Cartwright, D.I., Mantero, A.M., Steger, T.  and
Zappa, A., Groups acting simply transitively on the vertices of a
building of type~${\tilde A}_2$ I, II, {\it Geom. Ded.\/} {\bf 47} (1993),
143--166, 167--223.  \smallskip

\item{[\Ch]} Charney, R., Geodesic automation and growth functions
for Artin groups of finite type, Ohio State University preprint, (1993).
\smallskip

\item{[\CD]} Charney, R., and Davis, M., personal communication.
\smallskip

\item{[\Co]} Cohen, A., ``Recent results on Coxeter groups'',
preprint, (1993).

\item{[\D]} Dembowski, P., ``Finite Geometries", Ergebnisse der Mathematik
und ihrer Grenzgebiete, Band~{\bf44}, Springer-Verlag, Berlin,
Heidelberg, New York, 1968.
\smallskip

\item{[\Ep]} Epstein, D.B.A., Cannon, J.W., Holt, D.F., Levy, S.V.F.,
Paterson, M.S. and Thurston, W.P., ``Word processing in groups", Jones and 
Bartlett Publishers, Boston, 1992.   
\smallskip

\item{[\GS]} Gersten, S.M., and Short, H., Small cancellation theory and
automatic groups,  II, {\it Invent.\ Math.\/}, {\bf 105} (1991), 641--662.
\smallskip

\item{[\GSrat]}  Gersten, S.M., and Short, H., Rational subgroups of
biautomatic groups, Annals of Math. \ {\bf 134} 125 -- 158, 1991.

\item{[\Gr]} Gromov, M., Hyperbolic Groups in ``Essays in Group
Theory'', 75 -- 263, Gersten, S.M., ed., Springer Verlag,  M.S.R.I.
Series vol. 8, 1987.\smallskip

\item{[\Mou]} Moussong, G., ``Hyperbolic Coxeter Groups'', Doctoral
Thesis, Ohio State University, 1988.

\item{[\Pa]} Papasoglu, P., ``Geometric Methods in Group Theory'',
Doctoral Dissertation, Col\-umbia University, 1993.

\item{[\Pap]} Papasoglu, P., Strongly Geodesically Automatic Groups
Are Hyperbolic, University of Warwick preprint, 1994.

\item{[\Ro]} Ronan, M., ``Lectures on buildings'', Academic Press,
New York, 1989.

\item{[\T]} Tamaschke, O., ``Projektive Geometrie, I",
Biblio\-graph\-isches In\-sti\-tut, Mann\-heim, Wein, Z\"urich, 1969.

\bye